\documentclass[a4paper,12pt]{article}

\newtheorem{definition}{Definition}[section]
\newtheorem{proposition}{Proposition}[section]

\newtheorem{lemma}{Lemma}[section]
\newtheorem{theorem}{Theorem}[section]

\def\qed{\quad \vrule height7.5pt width4.17pt depth0pt} 
\usepackage{graphicx}
\usepackage{setspace}

\begin{document}

\doublespacing

\author{\large{Benjamin de Bivort}}
\title{Isotemporal Classes of $n$-gons}
\maketitle

Corresponding address:

Department of Molecular and Cellular Biology, Harvard University, 

7 Divinity Avenue, Cambridge, MA 02138

bivort@fas.harvard.edu, tel: 617-230-3769, fax: 617-495-8308

\pagebreak

\section{Abstract}
Here I present the present the first major result of a novel form of network analysis --- a temporal interpretation. Treating numerical edges labels as the time at which an interaction occurs between the two vertices comprising that edge generates a number of intriguing questions. For example, given the structure of a graph, how many ``fundamentally'' different temporally non-isomorphic forms are there, across all possible edge labelings. Specifically, two networks, $N$ and $M$, are considered to be in the same isotemporal class if there exists a function $\phi:N \to M$ that is a graph isomorphism and preserves all paths in $N$ with strictly increasing edge labels. I present a closed formula for the number of isotemporal classes ($\mathcal{N}(n)$) of $n$-gons. This result is strongly tied to number theoretic identities; in the case of $n$ odd, $\mathcal{N}(n)= \frac{1}{n} (\sum_{d \mid n}(2^{n/d -1}-1)\varphi(d)) )$, where $\varphi$ is the Euler totient function. 
\\
\\
\\
\\
\emph{Key Words}: temporal, network, graph, isomorphism, Euler, totient, phi, function, rotational, reflective, symmetry

\pagebreak

\section{Definitions and the Problem}

As a field, the mathematical analysis of networks has both sophistication and remarkable diversity. This is due largely to the surprising consistency with which a novel metaphorical interpretation of the values associated with the edges of a graph generates intriguing lines of mathematical investigation. To name a few examples, treating the values on edges as distances, throughput capacities, preferences, colors, etc. are metaphors for edge labels that have motivated profound advances. 

In this paper I present the first major result to come from a novel interpretation of labels. The values associated with an edge of a \emph{temporal network} indicate the time at which an interaction occurs between the two vertices comprising that edge. Many intriguing questions arise under this interpretation and here I present a classification of the number of fundamental temporal organizations possible on the $n$-gon class of graphs.

We'll begin with an example. Let graph theorists $A$, $B$, $C$, $D$ and $E$ belong to a strange academic society that meets often, but only two members at a time. Their most recent meeting history is given in Figure 1A. Forgetful Professor $B$ returned from a January trip to Paris bearing a miniature Eiffel Tower key-chain. Since then, he has misplaced the souvenir but clearly remembers lending it to another member of the society, although he is uncertain which. 

Precise conclusions based on the information provided in Figure 1A can be derived using a temporal interpretation of the network of interactions described therein. The precise definition of a temporal network:

\begin{definition}
A \textbf{temporal network} is a collection of vertices $V=\{v_1,$ $v_2,\ldots,$ $v_n\}$, edges $E=\{e_1,e_2,\ldots,e_m : \mbox{ where } e_i=\{v_j,v_k\} \mbox{ for some } v_j,v_k \in V\}$, \textbf{temporal labels} or \textbf{values} $T=\{t_1,t_2,\ldots,t_p : t_i \in \mathbf{R} \} \mbox{ where } n,m,p \in \mathbf{N}$, and a bijection, a \textbf{temporal labeling}, $\tau : E \to T$ from edges to temporal values.
\end{definition}

The temporal network of the interactions of $A$, $B$, $C$, $D$, and $E$ is shown in Figure 1B. If network 1B $=\{V,E,T,\tau\}$, then $V=\{A,B,C,D,E\}$, $E=\{ \{E,C\}, \{C,B\}, \{B,D\}, \{D,A\}, \{A,E\} \}$, and $T=\{$ April 1, April 9, June 15, July 8, July 10 $\}$. 

For convenience, we can order the times of the meetings (the temporal labels) in $T$, and replace the date of the first event with the number 1, the second date with the number 2 and so forth. This network is shown in Figure 1C. 

Desperate to recover his miniature, Prof. $B$ begins analyzing the temporal network. He concludes that any of his colleagues $D$, $E$ or $C$ could be in possession of the key-chain. $D$ could have it, had $B$ passed it to him on July 8 (time 4). Likewise $C$ could have it after the meeting at time 2; furthermore, $C$ could have passed it on to $E$ at their time 5 meeting. Only $A$ could not possibly possess the trinket, since $A$ met with $D$ and $E$ before either of those two could possibly have acquired it. 

That $E$ may possess the object reflects a temporal connectedness between $B$ and $E$ that motivates the following definition.

\begin{definition}
An ordered collection of vertices $P= \langle v_1, v_2, \ldots, v_m \rangle$ is a \textbf{temporal path} if $\tau(\{v_i,v_{i+1} \}) < \tau(\{v_j, v_{j+1} \} )$ for any $i < j$.
\end{definition}

After emailing $E$, $C$ and $D$ regarding the key-chain, the society's secretary contacted all the members to inform them that the April 1 meeting between $A$ and $D$ had been entered into the records incorrectly. The actual date was indeed in April, but which day in that month is unknown. 

If we assume that the ambiguous date was actually April 19, this changes the order of the interactions. The temporal network which corresponds to this alternative is shown in Figure 1D. However, inspection of these two networks reveals that every temporal path in network 1C is also a temporal path in network 1D. Indeed, the two networks are temporally isomorphic.

\begin{definition}
Let $N=\{V,E,T, \tau\}$ and $M=\{V', E', T', \tau'\}$ be temporal networks. If the function $\phi \colon V \to V'$ has the following properties then $\phi$ is a \textbf{temporal isomorphism} and $N$ and $M$ are \textbf{temporally isomorphic}:
\begin{itemize} 
\item $|V| = |V'|$
\item if $\{v_i, v_j\} \in E$, then $\{ \phi(v_i), \phi(v_j)\} \in E'$ (edge 	preservation)
\item if $P=\langle v_{a_1}, v_{a_2}, \ldots, v_{a_m} \rangle$ is a temporal 	path in N, then $\phi(P)$ $=\langle \phi(v_{a_1}),$ $\phi(v_{a_2}), 	\ldots,$ $\phi(v_{a_m}) \rangle$ is a temporal path in $M$ (temporal path 	preservation)
\end{itemize} 
\end{definition}

Any two networks which are temporally isomorphic are said to belong to the same \textbf{isotemporal class}, and if a particular function $\phi$ satisfies at least the edge preservation condition for networks $N$ and $M$, it is said to be a \textbf{graphical isomorphism} between the two networks. Temporal isomorphism and graphical isomorphism between $M$ and $N$ are denoted $M \cong_{\mathrm{T}} N$ and $M \cong_{\mathrm{G}} N$ respectively.

That two temporal networks (i.e. Figures 1C and 1D) can have fundamentally different temporal labelings, but belong to the same isotemporal class is an important property. Under a temporal interpretation, the temporal paths through a network (the paths over which an object could progress) are, in a sense, more fundamental descriptors of the network than the particular order in which the interactions occurred. 

An attempt to understand all the different temporal variants of a graph such as the 5-gon shown in Figure 1, would be well served by determining the number of 
5-gon different isotemporal classes. A more ambitious version of this question is, for a particular $n$, how many isotemporal classes ($\mathcal{N}$) of the $n$-gon are there?

\begin{definition}
A temporal network with vertices $V=\{v_1,v_2, \ldots v_n\}$ is an $n$-gon if for any $i \in \{1,2,\ldots, n\}$, $\{v_i,v_{i+1 \pmod{n}} \} \in E$, and $E$ contains no other edges.
\end{definition}

\section{Utility of the Line Graph}

\begin{definition}
Let $N=\{V,E,T,\tau\}$ be a temporal network with edges $E=\{e_1,e_2,\ldots,e_n\}$ . The \textbf{line graph of} $G$, $\mathcal{L}(G)$, is a graph with vertices $W=\{w_{e_1},w_{e_2},$ $\ldots,w_{e_n}\}$ corresponding to each edge in $G$. An edge exists in $\mathcal{L}(G)$ between vertices $w_{e_i}$ and $w_{e_j}$ if $e_i$ and $e_j$ share a vertex. The edge between $w_{e_i}$ and $w_{e_j}$ is directed toward $w_{e_i}$ if $\tau(e_i)>\tau(e_j)$, and toward $w_{e_j}$ if $\tau(e_i)<\tau(e_j)$. We write $w_{e_i} \rightarrow w_{e_j}$ if the edge between $w_{e_i}$ and $w_{e_j}$ is directed toward $w_{e_j}$
\end{definition}

The line graph of our example temporal network is shown in Figure 2A. Two line graphs are said to be \textbf{directionally isomorphic} ($\cong_{\mathrm{D}}$) if, in addition to edge preservation, there is preservation of the directedness of each edge. The line graph $\mathcal{L}(G)$ of a temporal $n$-gon provides a useful tool for counting the number of isotemporal classes because of the useful fact that for a temporal $n$-gon, $N$, $N \cong_{\mathrm{G}} \mathcal{L}(N)$. This follows immediately from the definition of the line graph; within any $n$-gon, one can inscribe another $n$-gon by rotation of $180/n$ degrees.

This fact is required to show that every isotemporal class of an $n$-gon can be uniquely and entirely described by a single directed line graph.

\begin{theorem}
Let $N$ and $M$ be temporal $n$-gons. $N \cong_{\mathrm{T}} M$ if and only if $\mathcal{L}(N) \cong_{\mathrm{D}} \mathcal{L}(M)$.
\end{theorem}
\emph{Proof} --- Let $N=\{V,E,T,\tau\}$, and $V=\{v_1,v_2, \ldots, v_n\}$. For the sake of simpler notation, whenever $i>n$ or $i<1$, $v_i$ shall be taken to mean $v_{i \pmod{n}}$ or $v_{i +nk \pmod{n}}$ respectively (where $nk+i>0$). First we show that $N \cong_{\mathrm{T}} M$ implies $\mathcal{L}(N) \cong_{\mathrm{D}} \mathcal{L}(M)$. For any $j$, either $P=\langle v_j, v_{j+1}, v_{j+2} \rangle$ or $Q=\langle v_{j+2}, v_{j+1}, v_{j} \rangle$ is a temporal path (as it must be the case that either $\tau(\{v_j, v_{j+1} \}) > \tau(\{v_{j+1}, v_{j+2} \})$ or vice-versa). Without loss of generality, we will assume that the former is a temporal path. By the definition of the line graph, $w_{ \{v_j, v_{j+1}\} } \rightarrow w_{ \{v_{j+1}, v_{j+2}\} }$. By assumption, there exists some $\phi \colon N \to M$ that is a temporal isomorphism, let $\phi(v_b)=u_b$ so that $\phi(P)=\langle u_{j}, u_{j+1}, u_{j+2} \rangle$. If the edge $\{u_c, u_{c+1}\}$ in $M$ corresponds to $x_{\{u_c, u_{c+1}\}}$ in $\mathcal{L}(M)$, then similarly, $x_{ \{u_j, u_{j+1}\} } \rightarrow x_{ \{u_{j+1}, u_{j+2}\} }$. Let $\psi \colon \mathcal{L}(N)\to \mathcal{L}(M)$ by $\psi(w_{\{v_i,v_{i+1}\}})=x_{\{ \phi(v_i), \phi(v_{i+1})\}}=x_{\{ u_i, u_{i+1}\}}$. This function clearly preserves edges, and since $x_{ \{u_j, u_{j+1}\} } \rightarrow x_{ \{u_{j+1}, u_{j+2}\} }$ whenever $w_{ \{v_j, v_{j+1}\} } \rightarrow w_{ \{v_{j+1}, v_{j+2}\} }$ it preserves directedness of the edges, and is a directional isomorphism.

To show the converse, that $\mathcal{L}(N) \cong_{\mathrm{D}} \mathcal{L}(M)$ implies $N \cong_{\mathrm{T}} M$, we invoke the fact that $N \cong_{\mathrm{G}} \mathcal{L}(N)$. Let the functions $\sigma_1 \colon N \to \mathcal{L}(N)$ by $\sigma_1(v_b)=w_{\{v_b,v_{b+1}\}}$, $\sigma_2 \colon M \to \mathcal{L}(M)$ by $\sigma_2(u_b)=x_{\{u_b,u_{b+1}\}}$, and $\psi \colon \mathcal{L}(N) \to \mathcal{L}(M)$ by $\psi(w_{\{v_c,v_{c+1}\}})=x_{\{u_d,u_{d+1}\}}$, be graphical isomorphisms. (We know that $\psi$ exists, by the assumption that $\mathcal{L}(N) \cong_{\mathrm{D}} \mathcal{L}(M)$). Let $\phi(v_i)=\sigma^{-1}_2(\psi(\sigma_1(v_i)))$. Graphical isomorphism is an equivalence relation, and so $\phi \colon N \to M$ will be a graphical isomorphism, since directional isomorphism implies graphical isomorphism. If $\phi$ preserves temporal paths, it will be a temporal isomorphism. Because $N$ is an $n$-gon, any temporal path will be in one of the following forms: $P=\langle v_c, v_{c+1}, \ldots, v_{c+m} \rangle$ or $\langle v_{c+m}, v_{c+m-1}, \ldots, v_{c} \rangle$. Without loss of generality, we will assume it is the former. By the definition of line graph, and application of $\sigma_1$, we know that $w_{\{v_c, v_{c+1}\}}$  $\rightarrow$ $w_{\{v_{c+1}, v_{c+2}\}}$ $\rightarrow \cdots$ $\rightarrow$ $w_{\{v_{c+m-1}, v_{c+m}\}}$, and since $\mathcal{L}(N) \cong_{\mathrm{D}} \mathcal{L}(M)$ by $\psi$, $x_{\{u_d, u_{d+1}\}}$  $\rightarrow$ $x_{\{u_{d+1}, u_{d+2}\}}$ $\rightarrow \cdots$ $\rightarrow$ $x_{\{u_{d+m-1}, u_{d+m}\}}$. The directedness of these edges implies that, after application of $\sigma^{-1}$, $\tau(\{u_d, u_{d+1}\}) <$ $\tau(\{u_{d+1}, u_{d+2}\}) < \ldots$  $< \tau(\{u_{d+m-1}, u_{d+m}\})$, where $\tau$ is the temporal labeling in $M$. Therefore, $\langle u_d, u_{d+1}, \ldots, u_{d+m} \rangle$ is a temporal path in $M$, and the function $\phi(v_c)=u_d$ is a temporal isomorphism from $N$ to $M$. 
\qed \\ 

This theorem places isotemporal classes in one-to-one and onto correspondence with isodirectional classes of line graphs. So, in order to determine $\mathcal{N}(n)$, we need only count the number of line graphs up to directional isomorphism. Given the trickiness of the counting arguments to come, we are well served to even further simplify our representation of isotemporal classes.

\begin{definition}
The \textbf{plus-minus form} ($\pm$-form) of an $n$-gon $N=$ $\{V,E,T,\tau\}$, is a $n$-gon labeled according to the following scheme. Noting the directedness of edges in $\mathcal{L}(N)$, edge $e_a$ receives a ``$+$'' label if $w_{e_{a-1}} \rightarrow w_{e_{a}}$ and $w_{e_{a}} \leftarrow w_{e_{a+1}}$, a ``$-$'' label if $w_{e_{a-1}} \leftarrow w_{e_{a}}$ and $w_{e_{a}} \rightarrow w_{e_{a+1}}$. If $w_{e_{a-1}} \rightarrow w_{e_{a}}$ and $w_{e_{a}} \rightarrow w_{e_{a+1}}$, or $w_{e_{a-1}} \leftarrow w_{e_{a}}$ and $w_{e_{a}} \leftarrow w_{e_{a+1}}$, then edge $e_a$ receives a ``0'' label. Thus, $\pm(N)=\{V,E,\{+,-,0\},f\}$, where $f \colon E \to \{+,-,0\}$ is the \textbf{$\pm$-labeling}.
\end{definition}

The $\pm$-form of our example temporal network is shown in Figure 2B.
There are several additional useful properties of the $\pm$-form of $n$-gons that follow directly from the definition.

\begin{itemize}
\item Let $A$ and $B$ be the $\pm$-form of temporal $n$-gons $N$ and $M$. 
$\phi \colon A \to B$ is a (\textbf{label}) isomorphism that preserves $\pm$ labels if and only if $\mathcal{L}(N) \cong_{\mathrm{D}} \mathcal{L}(M)$. 
\item There must be at least one edge of $A$ labeled with a $+$.
\item Any path through a $\pm$-form that starts and ends on edges labeled $+$, 
and containing no other $+$ labels, will have within it, precisely one edge labeled with a $-$.
\end{itemize}

\section{Let the Counting Begin}

Curiously this implies that in examining the labels of $\pm$-form in turn, we will find the $+$ and $-$ labels alternating, and interspaced by an arbitrary number of $0$ labels. You can see this pattern in Figure 2B. 

Here is our strategy for finding a formula for $\mathcal{N}(n)$, the number of isotemporal classes of an $n$-gon: 1) Count the number of distinct ways edges can be selected on an $n$-gon to receive non-zero labels. 2) Then, consider for each case, whether labeling an arbitrary first edge with a $+$ or $-$ label generates a different $\pm$-form. For the first part of this argument, we will need to invoke the help of the choose function. 

The number returned by the function ${{n} \choose {k}}$ can be interpreted as the number of order non-specific ways to select $k$ objects from a pool of $n$ distinct objects. If we let the pool of $n$ objects be the set $X=\{1,2,\ldots,n\}$, then ${{n} \choose {k}}$ returns the number of distinct subsets of $X$ of order $k$. Each of these subsets can be used to identify a class of labelings of a $\pm$-form of an $n$-gon by identifying those edges of the $n$-gon that are to receive non-zero labels. See Figure 3. The subset $\{3,4,6,8\}$ of $\{1,2, \ldots, 8\}$ represents those $\pm$-forms of the $8$-gon that have non-zero labels on the edges indicated in grey in Figure 3A. 

\begin{definition}
The \textbf{footprint} of the set $\{a_1, a_2, \ldots, a_k\}$ on an $n$-gon with edges $\{e_1, e_2, \ldots, e_n\}$ is the subgraph comprised of edges $\{e_{a_1}, e_{a_2}, \ldots,$ $e_{a_k} \}$.
\end{definition}

By no means does the choose function identify each distinct footprint uniquely, or even consistently. For example, all the footprints in Figure 3A are rotationally equivalent, and for this footprint, the choose function will identify 8 replicates. For the footprint shown in 3B only two replicates of the footprint will be identified. Additionally, the footprint in 3C is a mirror reflection of the first footprint of 3A; the two represent label isomorphic $\pm$-forms, but are identified by the choose function as distinct.

Four forms of symmetry will interfere with identifying distinct $\pm$-forms: mirror symmetry, skewed mirror symmetry, rotational symmetry, and skewed rotational symmetry. Examples, of $\pm$-forms and corresponding footprints of $n$-gons with eight of sixteen possible combinations of these types of symmetry are given in Figure 4A.

\begin{definition}
In an $n$-gon, a \textbf{vertex axis} of symmetry running through vertices $v_i$ and $v_{i+\frac{n}{2}}$, denoted $\overline{A}_{v_i}$ is an axis of \textbf{mirror symmetry} if $f(\{v_{i-a},$ $v_{i-a+1}\}) =$ $f(\{v_{i+a-1},v_{i+a}\})$. Similarly, $\overline{A}_{v_i}$ is an axis of \textbf{skewed mirror symmetry} if $f(\{v_{i-a},v_{i-a+1}\}) =$ $-f(\{v_{i+a-1},v_{i+a}\})$. 

An \textbf{edge axis} of symmetry running through $e_a$ and $e_{a+\frac{n}{2}}$ ($\overline{A}_{e_a}$) is an axis of mirror symmetry if $f(e_{a-k})=f(e_{a+k})$ and an axis of skewed mirror symmetry if $f(e_{a-k})=-f(e_{a+k})$.

A $n$-gon has \textbf{d-fold rotational symmetry} if for any edge $e_j$, $f(e_j)=f(e_{j+\frac{n}{d}})$ and \textbf{d-fold skewed rotational symmetry} if $f(e_j)=-f(e_{j+\frac{n}{d}})$.
\end{definition}

With those definitions, we can now approach the first task of our strategy, determining the number of distinct footprints:

\begin{theorem}
The number of footprints (up to reflective isomorphism) of an $n$-gon $N$, is  $\mathcal{M}(n)=\frac{1}{n} \displaystyle \sum_{d|n} (2^{n/d-1}-1)\varphi(d) $ if $n$ is odd, and $1+\frac{1}{n} \displaystyle \sum_{d|n} (2^{n/d-1}-1)\varphi(d)$ if $n$ is even. Here $\varphi(d)$ is the Euler totient function that returns the number of non-divisors of $d$.
\end{theorem}
\emph{Proof} --- We will begin with the odd case where $n=2k+1$. If each footprint indicates edges of $N$ that receive non-zero labels, it must contain $k=2, 4, \ldots, 2k$ edges, since the number of + labels must equal the number of - labels. Thus, the term $\sum_{i=1}^{k} {{n} \choose {2i}}$ counts all the footprints at least once. This can be simplified using basic binomial identities to $2^{n-1}-1$.

However, as we see in Figure 3, if a footprint lacks rotational symmetry it will be represented either $n$ or $2n$ times by choose, if it either lacks or has reflective symmetry respectively. And, if the footprint has at most $d$-fold rotational symmetry (as in Figure 3B), this term will identify it $n/d$ times. Since this formula does not claim to equate left and right-hand reflections of a footprint, we will only consider the mis-representation by choose of those footprints with rotational symmetry. 

It is our goal to compensate for the under-representation of rotationally symmetrical footprints by the choose function so that each footprint is counted either $n$ or $2n$ times depending on whether it has reflective symmetry. We will identify those footprints with at least $d$-fold symmetry with each term of the following formula: $\sum_{d \neq 1,d|n} \sum_{i=1}^{\frac{n/d-1}{2}} {{n/d} \choose {2i}} \Delta_d$. Again this simplifies to $\sum_{d \neq 1,d|n} (2^{n/d-1}-1) \Delta_d$. Here $\Delta_d$ is a correction factor specific to each $d$-fold symmetrical footprint that increases the number of occurrences of the under-represented class from $n/d$ to $n$.  When $\Delta_{d_1}$ corrects for each $d_1$-fold symmetrical footprint, the term also corrects to the same degree, all labelings with $d_1 d_2$-fold symmetry. 

Let us consider $p$-fold symmetrical labelings where $p_1$ is prime. As a prime, $p_1$ has no sub-divisors. Since, each application of the choose function will identify each $p_1$-fold symmetrical footprint $n/p_1$ times, and $n/p_1$ have already been identified by the initial $2^{n-1}-1$ term, $\Delta_{p_1}=p_1-1$, since $(p_1)\frac{n}{p}+\frac{n}{p}=n=\Delta_{p_1}(n/p) + (n/p)$. It is not a coincidence that $\Delta_{p_1}=\varphi(p_1)$.

We will prove that $\Delta_d=\varphi(d)$ by induction on the number of sub-divisors of $d$, and have already shown that when $d$ has no divisors, $\Delta_d=\varphi(d)$. So, assume that for all $d_i|d$, that $\Delta_{d_i}=\varphi(d_i)$. Since any $\Delta_{d_i}$ will contribute to the number of accumulated representations of footprints with $d$-fold symmetry, we can calculate $\Delta_d$ as follows: 
$$\frac{n}{d}\Delta_d = n-\frac{n}{d}(\displaystyle \sum_{d_i|d}\varphi(d_i) - \varphi(1) - \varphi(d) )$$
$$\Delta_d = d-\displaystyle \sum_{d_i|d}\varphi(d_i) + \varphi(1) + \varphi(d) -1$$ 
In these equations $\varphi(1)$ is subtracted since 1-fold symmetry corresponds to the rotationally asymmetrical case, which is accounted for by the $2^{n-1}-1$ term, and $\varphi(d)$ is subtracted since there is no previous term accounting for $d$-fold symmetry. Invoking the number theoretic fact that $n = \sum_{d|n} \varphi(d)$ to substitute and simplify, we have: $\Delta_d = d-d + 1 + \varphi(d) -1 = \varphi(d)$. So, if $\Delta_{d_i}= \varphi(d_i)$ for all divisors of $d$, then $\Delta_{d}= \varphi(d)$. This completes the second half of the proof by induction, and allows us, therefore, to combine terms for $\mathcal{M}(n)=$ $\frac{1}{n}((2^{n-1}-1)+\sum_{d\neq 1, d|n} (2^{n/d-1}-1)\varphi(d)) =$ $\frac{1}{n} \sum_{d|n} (2^{n/d-1}-1)\varphi(d)$, when $n$ is odd.

The proof of the even case of this formula is highly analogous, and for $n$ even, $\mathcal{M}(n)=$ $1+ \frac{1}{n} \sum_{d|n} (2^{n/d-1}-1)\varphi(d)$; the addition of 1 derives from the fact that for $n$ even, $\sum_{i=1}^{n/2} {{n} \choose {2i}}=2^{n-1}$.
\qed \\

What good is this formula, if it considers two footprints, isomorphic under reflection, to be distinct? As we will see, this result is sufficient to determine $\mathcal{N}(n)$ for $n$ odd. Furthermore, it is related to the number of binary necklaces fixed in the plane $\frac{1}{n} \sum_{d|n} \varphi(d)2^{n/d}$ [1]. Recall that our temporal networks are not ``fixed;'' labelings isomorphic under reflection are considered identical.

\begin{proposition}
Let $N$ be an $n$-gon with an axis of symmetry $\overline{A}_{e_a}$. $\overline{A}_{e_a}$ is an axis of mirror symmetry if and only if $f(e_a) \neq 0 \neq f(e_{a+\frac{n}{2}})$. $\overline{A}_{e_a}$ is an axis of skew symmetry if and only if $f(e_a) = f(e_{a+\frac{n}{2}})=0$. For vertex axes, $\overline{A}_{v_i}$ is an axis of symmetry if and only if it is an axis of skew symmetry
\end{proposition}
\emph{Proof} --- These properties follow directly from the fact that the $+$ and $-$ labels of $\pm(N)$, though potentially interspaced by 0 labels, must alternate in a +, -, +, - $\ldots$ fashion. 
\qed \\

Because every axis of symmetry in an \textbf{odd-gon} (an $n$-gon where $n=2k+1$) must pass through a vertex and an edge, it is impossible for an odd-gon to have an axis of mirror symmetry.  

\begin{proposition}
Let $N$ be a temporal $n$-gon, $\pm(N)=\{V,E,\{+,-,0\},f\}$ is label isomorphic to $-\pm(N)$ (the $\pm$- form with + and - labels interchanged) if and only if $\pm(N)$ has an axis of skewed mirror symmetry or d-fold skewed rotational symmetry. 
\end{proposition}
\emph{Proof} --- The set of automorphisms on an $n$-gon, the dihedral group $D_n$, consists of a flip, a rotation, or some combination of the two. Therefore if we let $\phi \colon \pm(N) \to -\pm(N)$ be a label isomorphism, then either $\phi(e_a)=e_{a+k}$ (a rotation), or $\phi(e_a)=e_{-a+k}$ (a flip and rotation). If $\phi$ satisfies the former, then, to preserve edge labels, $f(e_a)=-f(e_{a+k})$, and $\pm(N)$ has $n/k$-fold skewed rotational symmetry. Similarly, if $f(e_a)=-f(e_{-a+k})$, then  $\overline{A}_{e_{k/2}}$ is an axis of skewed mirror symmetry. 

Assume, without loss of generality, that $\pm(N)$ has an edge axis of skewed mirror symmetry;  $f(e_{b-k})=-f(e_{b+k})$ for some $b$ and any $k$. Let $\phi \colon \pm(N) \to -\pm(N)$ by $\phi(e_a)=e_{a+2(b-a)}$. Since these edges are symmetrically far from $e_b$, their labels will be + and -, or 0 and 0. Thus, $\phi$ will preserve edge labels from $\pm(N) \to -\pm(N)$.

If $\pm(N)$ has d-fold skewed rotational symmetry, then for any $j$, $f(e_j)=-f(e_{j+\frac{n}{d}})$. Let $\phi \colon \pm(N) \to -\pm(N)$ by $\phi(e_a)=e_{a+\frac{n}{d}}$. $\phi$ will clearly take + labels to - labels in $\pm(N)$, and + labels in $-\pm(N)$, and so constitutes a label isomorphism from $\pm(N)$ to $-\pm(N)$.
\qed \\

This proposition tells us exactly when alternatively labeling an arbitrary ``first edge'' of a footprint with a + or a - yields different $\pm$-forms: only when the $\pm$-form has neither skewed mirror symmetry nor skewed rotational symmetry. 

Further examination of the dihedral group and the choice of labeling the arbitrary first edge with a + or a - convinces us that the four cases of symmetry we have considered: skewed mirror, mirror, rotational and skewed rotational are indeed the only possible cases of symmetry that lead to miscounting by the choose function.

This lets us determine, for all combinations of symmetry, whether the choose function has mis-counted the number of isomorphically distinct footprints, and the number of distinct $\pm$-labelings (up to isomorphism: one or two) that each footprint needs to represent in our final formula (See Figure 5: Column A). 

Here, a ``1'' indicates a combination of symmetries such that the $\pm$-form of such a network, $P$, is directionally isomorphic to $-P$ ($P$ and $-P$ are identical); a ``2'' label indicates networks where $P$ is not isomorphic to $-P$ (and therefore, each footprint must represent 2 isotemporal classes). Figure 5: Column B gives the number of replicates of a particular footprint (again, up to isomorphism) identified by the formula of Theorem 2. Recall that left and right hand reflections were considered different footprints in that formula, so footprints without any reflective symmetry were counted twice.  

\begin{lemma}
Let $N$ be an odd-gon. $N$ cannot have $d$-fold skewed rotational symmetry.
\end{lemma}
\emph{Proof} --- Assuming the contrary, that there exists an $m$ such that for any edge $e_a$ in $\pm(N)$, $f(e_a)=-f(e_{a+m})$. Substituting $a+m$ for $a$ and so forth, we have, $f(e_a)=-f(e_{a+m})=f(e_{a+2m})=-f(e_{a+3m})=f(e_{a+4m})=,\cdots, =-f(e_{a+nm})$. But $e_{a+nm}=e_a$, so $f(e_a)=-f(e_{a+nm})=-f(e_a)$ is a contradiction.
\qed \\

\begin{theorem}

The number of isotemporal classes, $\mathcal{N}$, of an $n$-gon:\begin{itemize}
\item for $n=2k+1$, $\mathcal{N}(n)=\frac{1}{n} (\sum_{d \mid n}(2^{n/d -1}-1)\varphi(d)) )$

\item for $n=4k+2$, $\mathcal{N}(n)=\frac{1}{n} ( \sum_{d \mid n}2^{n/d -1}\varphi(d) - \sum_{c \mid \frac{n}{2}}2^{n/2c-1}\varphi(2c) ) + 2^{\frac{n-4}{2}} $

\item for $n=4k$, $\mathcal{N}(n)=\frac{1}{n} ( \sum_{d \mid n}2^{n/d -1}\varphi(d) - \sum_{c \mid \frac{n}{2}}2^{n/2c-1}\varphi(2c) ) + 2^{\frac{n-4}{2}} + 2^{\frac{n-8}{4}}- 2^{\lceil \frac{n-4}{8} \rceil -1}$
\end{itemize}
\end{theorem}
\emph{Proof} --- Let us first examine the $n$ odd case. By Lemma 3.1 and Proposition 3.1, we can eliminate any case of symmetry in which mirror or skewed rotational symmetry appear. Therefore only the first four rows of Figure 5 correspond to possible cases, and within these rows, the number of $\pm$-forms that correspond to a particular footprint (Column A) is identical to the number of copies of each footprint identified by the formula given in Theorem 2 (Column B). Therefore, that formula satisfies the odd case of this theorem. Needless to say, the even cases will be more complicated.

Since the odd-formula does not return the correct number (Column A) of $\pm$-labelings for four different categories of footprint (these are indicated with asterisks in Figure 5), additional correction terms are required. This correction will be done by adding or subtracting one replicate of each footprint in batches corresponding to cases of symmetry, so that after all the correction terms are taken into account, the sum of the counting terms of Columns B through F, across each row, will equal that in A.

 In Column C, for each $\pm$-form with mirror symmetry, another replicate is added. Column D subtracts a $\pm$-form replicate for each labeling with mirror and skewed mirror symmetries. Column E subtracts another $\pm$-form replicate for each labeling with skewed rotational symmetry, and finally Column F adds a $\pm$-form replicate for all labelings with skewed rotational and skewed mirror symmetries. The sum across each row of these correction terms and the initial value given by the odd-formula (Column B) is given in Column G. As Columns G and A are identical, implementing this sequence of corrections to the odd formula will yield the correct formula in the even cases; this is our road for the rest of the proof. \\

How are we going to count the number of footprints that have only skewed mirror symmetry and reflective skewed symmetries, or the number of footprints that have skewed rotational and rotational symmetries? In each case, a careful counting argument will give us the values we are interested in.
\\

\emph{Column C --- Adding a Replicate for Mirror Symmetry:}  Since an axis of mirror symmetry must pass through edges with non-zero labels, consider two polar edges of $N$ ``fixed.'' Each distinct half-footprint on one side of this axis or mirror symmetry will determine the footprint of the whole $n$-gon. Indeed, even the $\pm$ labels on the polar edges are determined by the $\pm$ labels on the rest of the half-footprint (if, say the first non-zero label away from the polar edge is a -, then that polar edge will be likewise flanked by a - on the other side, and therefore must have a + label itself). Each of the $\frac{n-2}{2}$ non-fixed edges in one half of $N$ may be independently included or not in a footprint, suggesting that there are $2^{\frac{n-2}{2}}$ possible footprints with mirror symmetry. However, if we let edges $e_1, e_2, \ldots, e_{(n-2)/2}$ be the candidate edges on one half of $N$, then the footprint $\{ e_{a_1}, e_{a_2}, \ldots, e_{a_k} \}$ will be isomorphic to $\{e_{(n-2)/2 + 1 -a_1}, e_{(n-2)/2 + 1 - a_2}, \ldots, e_{(n-2)/2 + 1 - a_k} \}$. That is, ``up'' and ``down'' oriented footprints are counted separately. Unfortunately, it is not sufficient to take $\frac{1}{2}2^{\frac{n-2}{2}}$ as the number of footprints, since it is possible that $\{ e_{a_1}, \ldots, e_{a_k} \} =$ $\{e_{(n-2)/2 + 1 -a_1}, \ldots, e_{(n-2)/2 + 1 - a_k} \}$ --- the case when a half-footprint has internal reflective symmetry. Such footprints are only identified once in the term $2^{\frac{n-2}{2}}$. If $n=4k+2$, this will occur $2^{\frac{n-2}{4}}$ times, since determining half of the half in question will determine the rest of the footprint. Similarly, if $n=4k$, the number of internally symmetrical half-footprints is $2^{\frac{n}{4}}$. Therefore if $n=4k$, the number of mirror symmetrical footprints is $\frac{1}{2} ( 2^{\frac{n-2}{2}} + 2^{\frac{n}{4}} ) =$ $2^{\frac{n-4}{2}} +2^{\frac{n-4}{4}}$, and if $n=4k+2$, the number of mirror-symmetrical footprints is $2^{\frac{n-4}{2}} +2^{\frac{n-6}{4}}$. These are our terms for Column C. \\

To count those footprints with both mirror and skewed mirror symmetry, we invoke a useful principle about multiple axes of symmetry in $n$-gons.  \\

\emph{Mirror and Skewed Mirror Axes Reflect Each Other:} To avoid repetition, we will address edge skewed mirror axes; proofs of these properties for vertex axes of skewed mirror symmetry are analogous. Let $\overline{A}_{e_a}$ be a skewed mirror axis, and $\overline{A}_{e_{a-k}}$ be a mirror axis. By definition of a mirror axis, $f(e_{a-k-b})=f(e_{a-k+b})$. Reflecting these edges through the skewed axis of symmetry we find that $f(e_{a-k-b})=-f(e_{a+k+b})=-f(e_{a+k-b})$. Therefore, $\overline{A}_{e_{a+k}}$ is an axis of mirror symmetry. In this way mirror and skewed mirror axes reflect each other. Consequently, if there are $l$ axes of mirror symmetry in $N$ and at least one axis of skewed mirror symmetry, then there will be $l$ axes of skewed mirror symmetry, each found halfway between two adjacent axes of mirror symmetry. This implies in turn that there if there are both mirror and skewed mirror axes in $N$, there must be, in total, even number of reflective axes, two of which are perpendicular.

\emph{Columns D and F Cancel Each Other Out:} Column C requires counting the number of footprints that correspond to $\pm$-forms with axes of both mirror and skewed mirror symmetry. Since axes of mirror and skewed mirror symmetry alternate, we can consider the subgraph of $N$ between a mirror axis and the nearest skewed mirror axis moving in a clockwise direction. Let this subgraph have footprint $A=\langle a_1,a_2,\ldots, a_g \rangle$ and $\pm$-labels $f(A)$. Let the reflection of $A$, be $A'=\langle a_g, \ldots, a_2, a_1 \rangle$. Proceeding clockwise around $N$ to the next subgraph between two axes, $A$ is reflected and negated across the skewed axis to yield $-f(A')$. Similarly $-f(A')$ is reflected around the next mirror axis to give $-f(A)$, and so forth. Thus $N$ is comprised of adjacent subgraphs with labels $f(A), -f(A'), -f(A), f(A'), f(A), -f(A'), \ldots$. This graph has $n/2r$-fold skewed rotational symmetry, where $r$ is the number of edges between adjacent mirror and skewed mirror axes, and rotating $N$ by $2r$ edges takes $f(A)$ to $-f(A)$ and $-f(A')$ to $f(A')$. 

Likewise, if $N$ has d-fold skewed rotational symmetry and skewed mirror symmetry, it will also have an axis of mirror symmetry. Let this axis of skewed mirror symmetry $\overline{A}_{e_a}$ partition $N$ into halves $f(A)$ and $-f(A')$.  By definition of skewed rotational symmetry, $f(B)=\langle f(e_{a-n/d}),$ $f(e_{a-n/d+1}), \ldots,$ $f(e_a) \rangle =$ $\langle -f(e_{a+1}), -f(e_{a+2}), \ldots,$ $-f(e_{a+n/d}) \rangle =-f(B)$. Assuming $B \subset A$, then $-f(B)=-f(B')$, and $f(B)=f(B')$, or $f(B)$ has internal mirror symmetry. 

Therefore, $N$ has skewed mirror and mirror symmetries if and only if it has 
skewed mirror and skewed rotational symmetries. Thus, the cases to be identified in Columns D and F are one in the same, and $\pm$-forms with only skewed mirror and mirror symmetries, and likewise $\pm$-forms with only skewed rotational and skewed mirror symmetry cannot exist, since they both imply the existence of the third type of symmetry. These cases are marked by double asterisks in Figure 5. Since in Column D we were to subtract the number of such cases, while adding them in Column F, the net contribution of the correction terms generated by these two columns is zero. \\

This property is remarkably convenient. All we need now is the number of $\pm$-forms with skew rotational symmetry.
\\

\emph{Column E --- Subtracting a Replicate for Skewed Rotation:} $N$ has $d$-fold skewed rotational symmetry if and only if the footprint of $N$ has $d$-fold rotational symmetry, $d$ is even, and the number of non-zero labels on $n/d$ adjacent edges is odd. This property is self-apparent, when it is considered that there must be an even number of non-zero labels, and that if $n/d$ adjacent edges were to contain an even number of edges with non-zero labels, the labels on those edges would then be $C=\langle 0_1, 0_2, \ldots, +, \ldots, -, \ldots, 0_{n/d-1}, 0_{n/d} \rangle$, a sequence identical to that found on the next $n/d$ edges. If $C$ contains an odd number of non-zeros, $C=\langle 0_1, 0_2, \ldots, +, \ldots, +, \ldots, 0_{n/d-1}, 0_{n/d} \rangle$, then the sequence of labels on the next $n/d$ adjacent edges would be $-C=\langle 0_1, 0_2, \ldots, -, \ldots, -, \ldots, 0_{n/d-1}, 0_{n/d} \rangle$, satisfying the definition of skewed rotational symmetry.

In order to count the number of skewed rotational footprints, we will need to use a similar argument as that used in Theorem 2. Summing over possible even $c$-folds, the number of ways to select an odd number of edges from $n/c$ edges is $\sum_{c \mid \frac{n}{2}} \sum_{k=0}^{\lceil \frac{n}{4c} - 1 \rceil} {{\frac{n}{2c}} \choose {2k+1}} = \sum_{c \mid \frac{n}{2}} 2^{n/2c-1}$. In order to count each occurrence $n/2$ times, we must introduce a correction factor similar to $\Delta_d$. An argument analogous to that given in the proof of Theorem 2 shows that $\frac{2}{n} \sum_{c \mid \frac{n}{2}} 2^{n/2c-1} \varphi(2c)$ returns the number of $c$-fold skewed rotationally symmetrical $n$-gons, ignoring, as was ignored in Theorem 2, the double counting of footprints that lack axes of reflective symmetry. 

Column E demands that we subtract from our formula only one replicate of each $\pm$-form that has skewed rotational symmetry. Therefore the term $\frac{1}{2} ( (\frac{2}{n} \sum_{c \mid \frac{n}{2}} 2^{n/2c-1} \varphi(2c)) + \Lambda )$ will give the number of cases taking into account reflective asymmetry. Here $\Lambda$ is the number of $\pm$-forms with skewed rotational symmetry, and some kind of reflective symmetry (i.e., those that are only counted once by the summation term). As we saw above, if $\pm$-form has skewed rotational symmetry and some form of reflective symmetry, then it must have both mirror and skewed mirror symmetries. And with both kinds of reflective symmetry present, $N$ must contain at least two perpendicular axes of symmetry.  

In general, we must consider the possibility that the axis perpendicular to the mirror axis could be either another axis of mirror symmetry, or an axis of skew symmetry, and that these two cases need to be counted separately. Let the number of $\pm$-forms with at least two axes of symmetry (our correction factor) be $\Lambda = \Lambda_{\mathrm{skew}} + \Lambda_{\mathrm{mirror}}$, the sum of the number of $\pm$-forms where the perpendicular axis is a skewed mirror axis or a mirror axis, respectively. 

If $n=4k+2$, the axis perpendicular to the axis of mirror symmetry must be a skewed mirror axis, as it passes through vertices. Here, $\Lambda_{\mathrm{mirror}}=0$, and since all other cases have a perpendicular skewed mirror axis, determining the number of quarter-footprints will determine the number of $\pm$-gons with this form ($\Lambda_{\mathrm{skew}})$. There are $\frac{n-2}{4}$ edges in this quadrant which can be independently included or not in the quarter-footprint. Therefore $\Lambda_{\mathrm{skew}}=2^{\frac{n-2}{4}}=\Lambda$ is the number of $n=4k+2$-gons with skewed rotational symmetry and reflective symmetry.

If $n=4k$, the axis perpendicular to the mirror axis placed by assumption can be either an edge skewed mirror axis, or an edge mirror axis. Assuming the former, $\Lambda_{\mathrm{skew}}=2^{\frac{n-4}{4}}$, as there are $\frac{n-4}{4}$ edges in the quadrant between the mirror axis, and the edge skewed mirror axis. 

If the perpendicular axis is another mirror axis, then, because skewed mirror and mirror axes must alternate, there must be an odd number of axes between the perpendicular mirror axes.  Therefore, one axis (either mirror or skewed mirror) must be a bisecting axis between the two perpendicular axes of mirror symmetry. With the existence of this bisecting axis guaranteed, determining one octant of the footprint of $N$ determines the whole, and $\Lambda_{\mathrm{mirror}}=2^{\lceil \frac{n-4}{8} \rceil}$. Here the ceiling function accounts for the two possibilities: $n=8k$ or $n=8k+4$, as there are $\lceil \frac{n-4}{8} \rceil$ edges in one octant of $n$-gons of both of these types. Thus, $\Lambda = 2^{\frac{n-4}{4}} + 2^{\lceil \frac{n-4}{8} \rceil}$.

Therefore, the number of $\pm$-forms with skewed rotational symmetry is $\frac{1}{2} ( (\frac{2}{n} \sum_{c \mid \frac{n}{2}} 2^{n/2c-1} \varphi(2c)) + 2^{\frac{n-2}{4}})$ if $n=4k+2$. If $n=4k$, then the number of skewed rotational cases is $\frac{1}{2} ( (\frac{2}{n} \sum_{c \mid \frac{n}{2}} 2^{n/2c-1} \varphi(2c)) + 2^{\frac{n-4}{4}} + 2^{\lceil \frac{n-4}{8} \rceil})$.

\emph{Piecing It All Together:} Assembling terms from the odd-formula, and Columns C and E, for $n=4k+2$, we have $\mathcal{N}(n)=$
$$\frac{1}{n}\sum_{d \mid n} (2^{n/d-1}-1)\varphi(d) +1 +(2^{\frac{n-4}{2}} + 2^{\frac{n-6}{4}}) - \frac{1}{2}( \frac{2}{n}(\sum_{c \mid \frac{n}{2}} 2^{n/2c-1}\varphi(2c))+2^{\frac{n-2}{4}} )$$
$$\mathcal{N}(n)= \frac{1}{n} (\sum_{d \mid n} 2^{n/d-1}\varphi(d) - \sum_{c \mid \frac{n}{2}}2^{n/2c-1}\varphi(2c) ) + 2^{\frac{n-4}{2}}$$

For $n=4k$, (slightly more complicated of course), we have $\mathcal{N}(n)=$
$$\frac{1}{n}\sum_{d \mid n} (2^{n/d-1}-1)\varphi(d) +1 +(2^{\frac{n-4}{2}} + 2^{\frac{n-4}{4}}) - \frac{1}{2}( \frac{2}{n}(\sum_{c \mid \frac{n}{2}} 2^{n/2c-1}\varphi(2c))+2^{\frac{n-4}{4}} + 2^{\lceil \frac{n-4}{8} \rceil} )$$
$$\mathcal{N}(n)= \frac{1}{n} (\sum_{d \mid n} 2^{n/d-1}\varphi(d) - \sum_{c \mid \frac{n}{2}}2^{n/2c-1}\varphi(2c) ) + 2^{\frac{n-4}{2}}  + 2^{\frac{n-8}{4}} - 2^{\lceil \frac{n-4}{8} \rceil -1}$$
\qed\\

The first 25 terms of the sequence defined by this result for $n=3, 4, \ldots 27$ are: 1, 3, 3, 8, 9, 20, 29, 60, 93, 188, 315, 618, 1095, 2118, 3855, 7414, 13797, 26482, 49939, 95838, 182361, 350580, 671091, 1292604, 2485533.

There are several intriguing lines of investigation suggested by this result, such as the relationship between temporal networks and $k$-ary necklaces, and colored graphs more generally. This proof also provides a new demonstration of the number theoretic result that $\sum_{d|n} \varphi(d)2^{n/d}$ is divisible by $n$. Furthermore, the above sequence of integers appears to be converging to 2; Showing that $\lim_{n \rightarrow \infty} \mathcal{N}(n) / \mathcal{N}(n-1) = 2$ would reveal a surprising result: adding a new member to a temporal cycle of interactions only doubles the number of fundamental ways objects can flow through the system, even while the number of distinct interaction sequences grows roughly by a factor of $n$. 

The novel field of temporal networks is ripe with open and undiscovered questions, such as the number of isotemporal classes of other graph forms, characterization of temporal path properties such as Hamiltonian restrictions, and algebraic implications of temporal label permutations that preserve temporal isomorphism. 
\\
\\
\\
\\
\emph{Reference}
\\
1. N. J. Fine, Classes of periodic sequences, \emph{Illinois J. Math}., 2 (1958), 285-302.
\\
\\
\small{\emph{Acknowledgments}. This paper would not have been possible without the continued guidance and insightful comments of Tristan Tager and Professor David Kraines of Duke University, and for this help I am very grateful.}
\\
\\
\\
\\
\\
\begin{figure}[h]\begin{center}\includegraphics {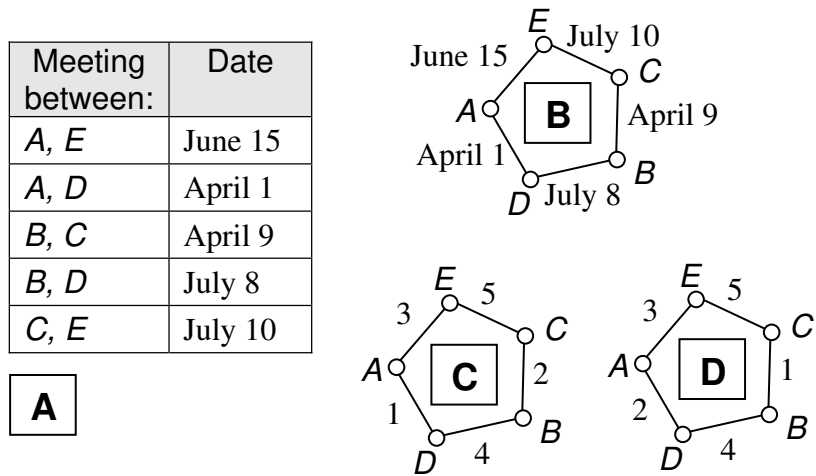} \end{center}
\caption{ \small{\emph{Example of a temporal network} --- \textbf{A}: Members of a curious academic society, Profs. $A$, $B$, $C$, $D$ and $E$ only meet in pairs, on the dates shown. \textbf{B}: Temporal network representing the interactions described in the table with date edge labels. The network in \textbf{C} replaces the dates with the ordered rank of the events, and \textbf{D} gives the alternative interpretation of the time of the interaction between A and D (See text).} }
\end{figure}

\begin{figure}[h]
\begin{center}\includegraphics{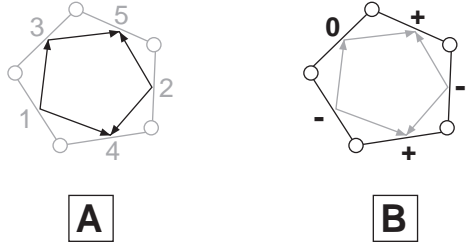} \end{center}\caption{ \small{\emph{Line Graph of the Temporal Networks from Figure 1} --- \textbf{A}: The temporal network from Figure 1C is shown in grey beneath the line graph of that network. \textbf{B}: The $\pm$-form of that line graph is shown in black above the line graph.} }
\end{figure}

\begin{figure}[h]
\begin{center}\includegraphics{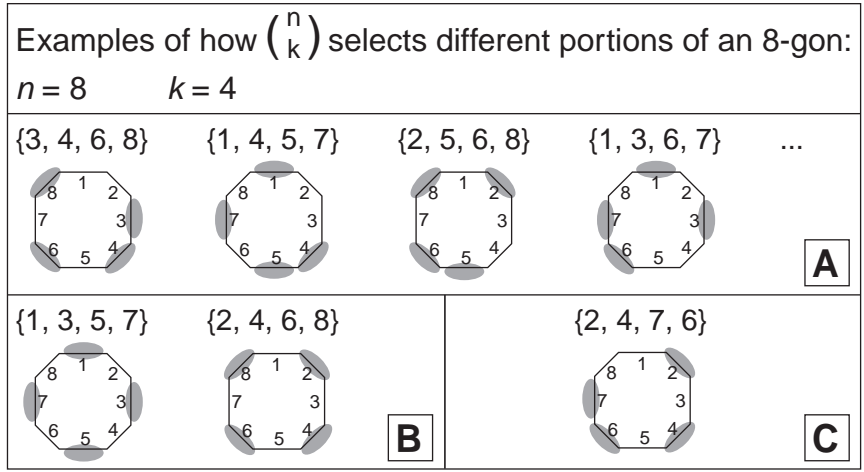} \end{center}
\caption{ \small{\emph{Examples of Footprints on an 8-gon} --- \textbf{A}: Rotational isomorphs of footprints without reflective symmetry are counted $n$ times in the ${{n} \choose {k}}$ term that identifies them. \textbf{B}: Examples of footprints of the 8-gon with 4-fold rotational symmetry. Footprints with $d$-fold rotational symmetry are represented $n/d$ times. \textbf{C}: Those footprints without reflective symmetry are additionally counted for their ``left'' and ``right-hand'' versions. The footprint in \textbf{C} is isomorphic to all those in \textbf{A}, but is identified by a distinct footprint. } }
\end{figure}

\begin{figure}[h]
\begin{center}\includegraphics{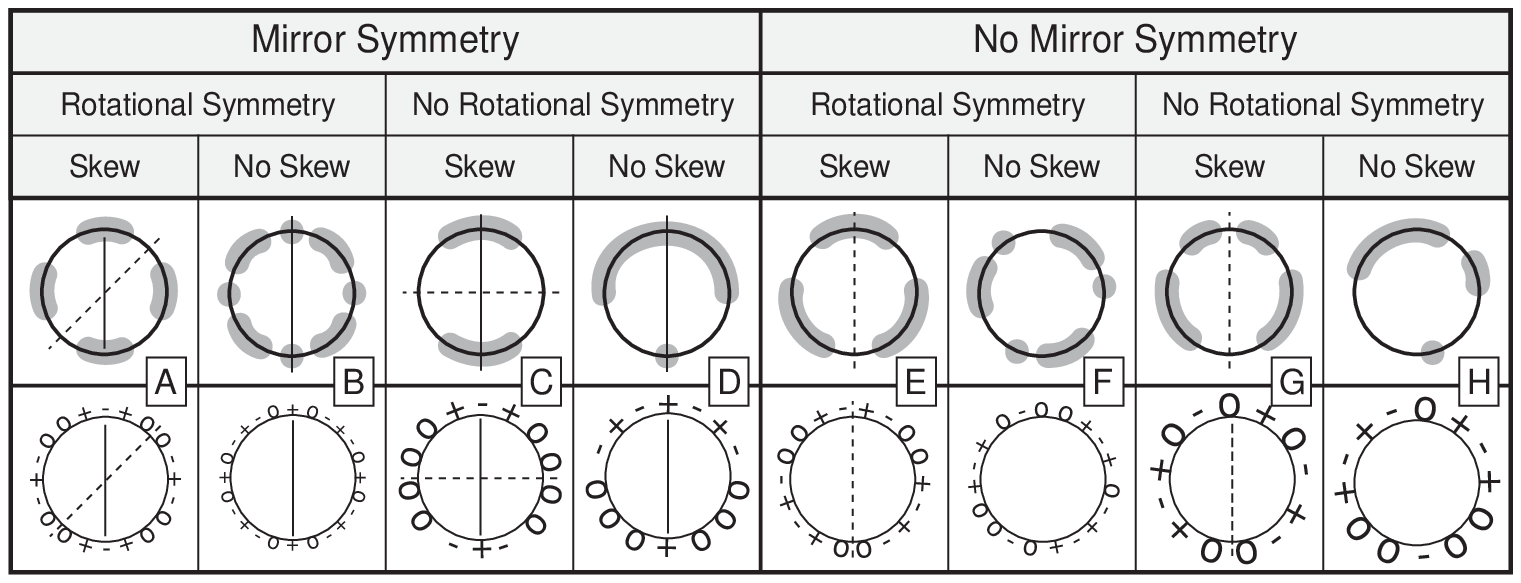} \end{center}
\caption{ \small{\emph{Examples of Various Combinations of Symmetry} --- 
Shown for eight of sixteen possible combinations of symmetry, is an example of an $n$-gon with a particular type of symmetry, with its corresponding footprint above. Dashed lines indicate axes of mirror symmetry, and solid lines axes of skewed mirror symmetry. } }
\end{figure}

\begin{figure}[h]
\begin{center}\includegraphics{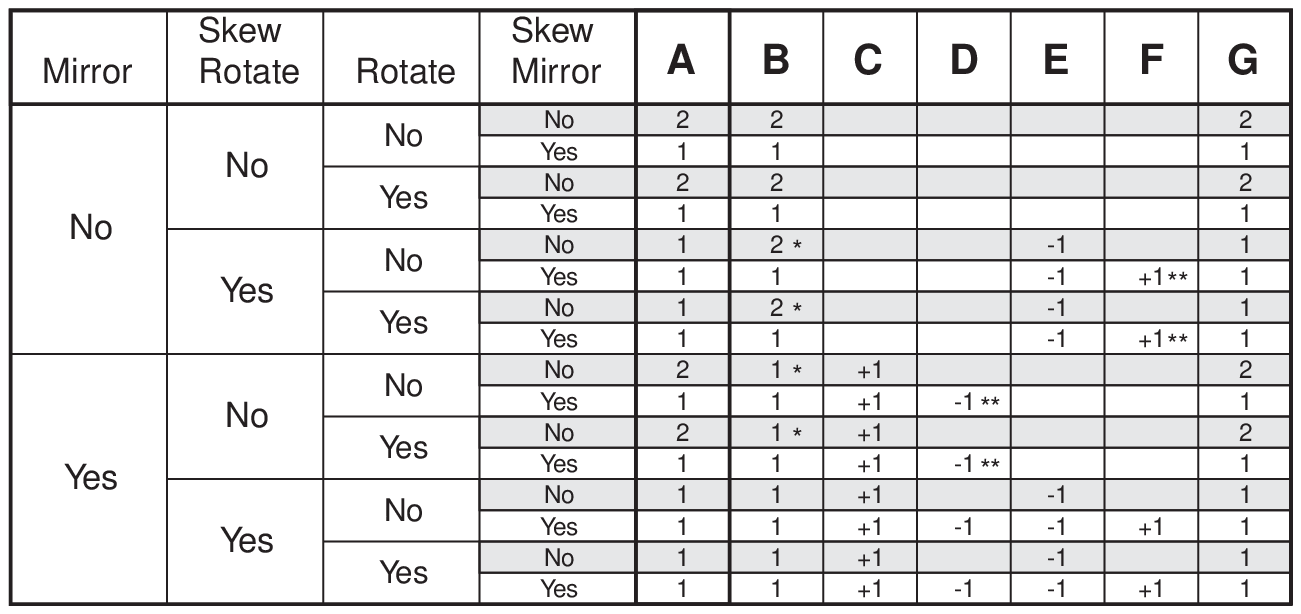} \end{center}
\caption{ \small{\emph{All Possible Symmetries and a Strategy for Calculating the Even Formula} --- Column \textbf{A} gives the number of isotemporal classes each footprint must ultimately represent. The number of replicates of a particular footprint, as identified by the odd-formula is given in \textbf{B}. \textbf{C}, \textbf{D}, \textbf{E}, and \textbf{F} represent corrections taken to revise the values in \textbf{B} to equal those in \textbf{A}: respectively, addition of mirror symmetric cases, subtraction of mirror and skewed mirror symmetric cases, subtraction of skewed rotational cases, and addition of skewed mirror and skewed rotational cases. \textbf{G} gives the sum of \textbf{B} through \textbf{F}, and as the correction strategy is sound, has entries equal to the goal of \textbf{A}. $\ast$: rows where the odd-formula and the goal differ. $\ast \ast$: combinations of symmetry that are impossible. } }
\end{figure}

\end{document}